\documentclass[twoside]{ae100prg}
\usepackage{booktabs}

\usepackage[pdftex]{graphicx}
\usepackage{color}
\usepackage{amstext,amssymb,amsfonts,amsmath,amsbsy,latexsym}
\newcommand{\Ric}{\operatorname{Ric}}
\newcommand{\Scal}{\operatorname{Scal}}

\begin{document}

\title{Geometrostatics: the geometry of static space-times}

\author{Carla Cederbaum}

\address{Mathematics Department, Duke University, Box 90320, Durham,
NC 27708-0320, USA}

\email{carla@math.duke.edu}

\begin{abstract}
We present a new geometric approach to the study of static isolated general
relativistic systems for which we suggest the name \emph{geometrostatics}. After
describing the setup, we introduce localized formulas for the ADM-mass and
ADM/CMC-center of mass of geometrostatic systems. We then explain the
pseudo-Newtonian character of these formulas and show that they converge to
Newtonian mass and center of mass in the Newtonian limit, respectively, using
Ehlers' frame theory. Moreover, we present a novel physical interpretation of
the level sets of the canonical lapse function and apply it to prove uniqueness
results. Finally, we suggest a notion of force on test particles in
geometrostatic space-times.
\end{abstract}

\section{Introduction}
Static isolated general relativistic systems have been studied from a number of
perspectives including their regularity, compactification and
asymptotic considerations, symmetry classifications, construction of explicit
solutions etc. They serve as models of static stars and black
holes. Also, they play an important role in R.~Bartnik's definition of
quasi-local mass and his associated conjecture on static metric extensions
\cite{Szabados}.

\begin{figure}[h]
\begin{center}
\includegraphics[scale=0.15,angle=-90]{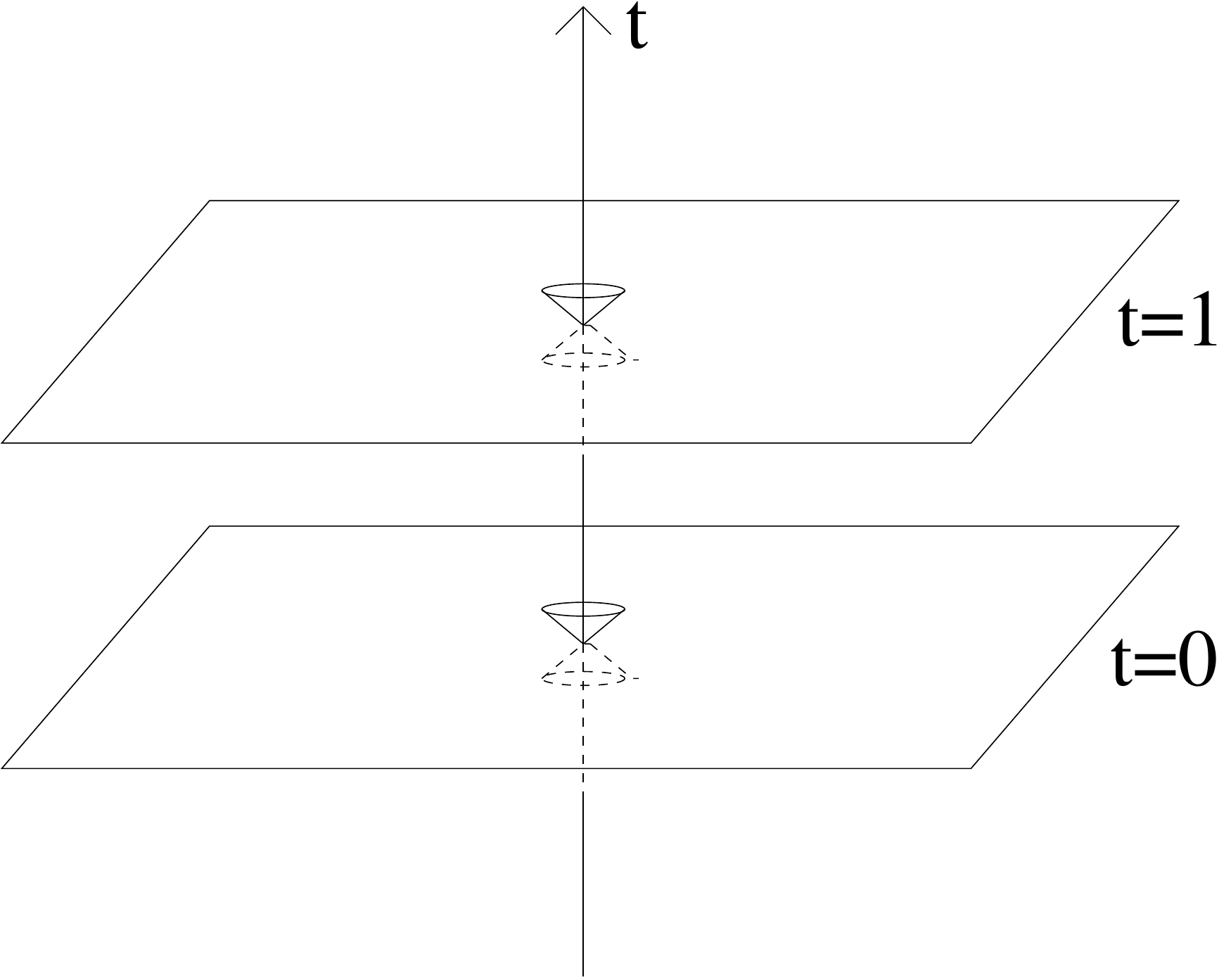}
\caption{\label{timeslice}The time-slices of a canonically decomposed static space-time.}
\end{center}
\end{figure}

Here, we present a new geometric approach to the study of static isolated
systems and their physical properties for which we suggest the
name \emph{geometrostatics}. We consider space-times that are
\emph{static} (possess a smooth global time-like Killing vector
field that is hypersurface-orthogonal) and \emph{isolated} (see below). Static
space-times generically possess a $3+1$-decomposition with vanishing shift
vector. In this \emph{canonical} decomposition, the \emph{canonical} lapse function is
time-independent and coincides with the Lorentzian length of the time-like
Killing vector field. The spacelike time-slices orthogonal to the
time-like Killing vector field are all isometric and have vanishing extrinsic curvature, see Figure \ref{timeslice}. Their induced Riemannian metric is time-independent. We will subsequently identify all canonical time-slices.

\vspace{-2ex}

For our purposes, static systems are called \emph{isolated} if the Riemannian
metric and the lapse function on the time-slice decay suitably fast
to the flat metric and the constant $1$, respectively, at spacelike infinity see \cite{Cederbaum} for a precise definition of our asymptotic flatness condition in the language of weighted Sobolev spaces. Moreover, we request that the space-time satisfies the \emph{vacuum} Einstein equations
outside some spatially compact tube in the space-time (or, in other words,
outside some compact set in the time-slice). This can be interpreted as a
(spatially) finite extension of the sources, whether they are matter sources
and/or black holes.

This article is structured as follows: In Section \ref{setup}, we will
introduce the central equations of geometrostatics and summarize a few of their
central analytic properties. In Section
\ref{levels}, we present a novel physical interpretation of the level sets of
the lapse function of a geometrostatic system and discuss some applications of
this insight. In Section \ref{pseudo}, we will perform a conformal
transformation into what we suggest to call \emph{pseudo-Newtonian} variables.
Moreover, we will define and analyze localized surface integral expressions
for the mass and center of mass of geometrostatic systems. Finally, in Section
\ref{limit}, we will discuss the Newtonian limit of geometrostatics.

Further details can be found in my thesis \cite{Cederbaum}. 

\section{Geometrostatics}\label{setup}
The Lorentzian metric of a generic static space-time $\mathbb{R}\times M^3$ can globally be decomposed as
\begin{equation}
ds^2=-N^2c^2dt^2+g,
\end{equation}
where $N:=\sqrt{-ds^2(\partial_{t},\partial_{t})}$ is the (canonical) lapse function arising as the Lorentzian length of the time-like Killing vector field $\partial_{t}$, $c$ is the speed of light, and $g$ is the Riemannian metric induced on the time-slice $M^3$. Observe that $N$ is non-negative and vanishes only along Killing horizons.

In the vacuum region outside the matter, the (vacuum) Einstein equations imply that these variables satisfy the
so-called \textit{vacuum static metric equations}
\begin{eqnarray}\label{SME1}
 N\Ric&=&\nabla^2N\\\label{SME2}
\triangle N&=&0,
\end{eqnarray}
where $\Ric$ is the Ricci curvature tensor of $g$, $\nabla^2N$ denotes the covariant
Hessian, and $\triangle N$ denotes the covariant Laplacian of $N$ with respect to $g$. It is well-known that solutions to the vacuum static metric equations are real analytic in suitable coordinates \cite{MzH}. 

We define a \emph{geometrostatic system} to be an asymptotically flat Riemannian $3$-manifold $(M^3,g)$ endowed with a smooth positive lapse function $N$ so that the vacuum static metric equations \eqref{SME1} and \eqref{SME2} are satisfied. Hence, geometrostatic systems model the vacuum region outside the support of the matter and the horizons of all black holes  within a slice of an asymptotically flat static space-time. The lapse function $N$ describes the lapse of time in the space-time.

\section{The Level Sets of the Lapse}\label{levels}
\begin{figure}[h]
\begin{center}
\resizebox{35ex}{!}{\input{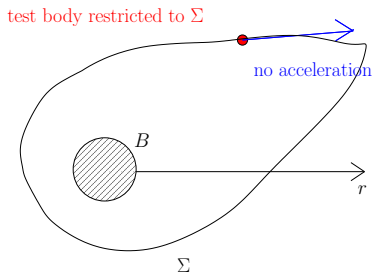}}
\caption{\label{levelset}A test particle constrained to a surface $\Sigma$.}
\end{center}
\end{figure}

In Newtonian gravity, the relevant gravitational variable is the Newtonian potential. The gradient of the potential defines the force on a unit mass test particle. This has a well-known consequence for the \emph{equipotential surfaces} (or level sets of the Newtonian potential): if a test particle is constrained into one of these surfaces then the gravitational force does not have a tangential component and hence the test particle does not tangentially accelerate within the surface, see Figure \ref{levelset}.

Surprisingly, the ``same'' is true for level sets of the lapse function in a geometrostatic system. In order to make this rigorous, we make the following definitions: Consider a closed smooth surface $\Sigma\subset M^3$ in a geometrostatic system $(M^3,g,N)$ arising as the $n=n_{0}$ level set of a smooth function $n:M^3\to\mathbb{R}$. A time-like\footnote{Here, time-like curves and the time functional are taken with respect to the static space-time metric $ds^2=-N^2c^2dt^2+g$ induced by $(M^3,g,N)$.} curve 
\begin{equation}
\mu(\tau)=(t(\tau),x(\tau))
\end{equation}
satisfying $x(\tau)\in\Sigma$ is called a \emph{test particle constrained to $\Sigma$} if it is a critical point of the time functional 
\begin{equation}
T(\mu):=\int_{\tau_{0}}^{\tau_{1}} \lbrace\vert\dot{\mu}(\tau)\vert+\sigma(n\circ x(\tau)-n_0)\rbrace\,d\tau,
\end{equation}
where $\sigma\in\mathbb{R}$ is a Lagrange multiplier ensuring that all comparison curves are also constrained to $\Sigma$. With this notion of constrained test particle, we say that a smooth closed surface $\Sigma$ is an \emph{equipotential surface}  if every test particle constrained to $\Sigma$ is a geodesic in $\Sigma$ with respect to the induced $2$-metric, see Figure \ref{levelset}. Analyzing the geodesic equation, we find that a surface $\Sigma\subset M^3$ is an equipotential surface in $(M^3,g,N)$ if and only if $\Sigma$ is a level set of $N$. Thus, the level sets of the lapse function $N$ in geometrostatics play precisely the same role as those of the Newtonian potential in Newtonian gravity.

In a static vacuum space-time, the Einstein constraint equations reduce to $\Scal=0$. In particular, the lapse function does not appear in this constraint equation. As a consequence, Choquet-Bruhat's theorem (see e.g.~\cite{CHB}) on the local existence and uniqueness of solutions to the Einstein equations implies that the space-time induced by $(M^3,g,N)$ is in fact \emph{independent}\footnote{This assumes that the lapse function \emph{exists} in the first place.} of the lapse function $N$. Combining this view of the lapse function with the physical interpretation of its level sets of the lapse as well as with the vacuum static metric equations \eqref{SME1} and \eqref{SME2} and the assumed asymptotic conditions for $g$ and $N$, we obtain that the lapse function is indeed unique if it exists. This is to say that if $(M^3,g,N)$ and $(M^3,g,\widetilde{N})$ are geometrostatic systems, then $N=\widetilde{N}$.
We interpret this result as saying that ``there is only one way of synchronizing time at different locations in a geometrostatic space-time such that one sees staticity'' just as, for a geodesic, ``there is only one way of walking along a geodesic such that one does not accelerate (up to affine transformations of the curve parameter)''. The affine freedom of the parameter along the geodesic does not make an appearance in the geometrostatic space-time picture because we fixed the lapse function to asymptotically converge to $1$ at spacelike infinity and therewith fixed the time unit.  

\section{Pseudo-Newtonian Gravity}\label{pseudo}
The geometrostatic variables $g$ and $N$ are ideal for investigating geometric and relativistic effects influencing test particle behavior and the behavior of light rays. In order to better understand asymptotic and analytic properties of solutions, however, it is more convenient to perform a conformal change and consider the new variables 
\begin{eqnarray}\label{gamma}
\gamma&:=&N^2g\\\label{U}
U&:=&c^2\ln N.
\end{eqnarray}
These variables have been used by many authors\footnote{albeit without explicit reference to the speed of light.}, see e.g.~\cite{KM}. We suggest to call them \emph{pseudo-Newtonian metric and potential}, respectively. The vacuum static metric equations \eqref{SME1}, \eqref{SME2} translate into
\begin{eqnarray}\label{conf1}
\Ric_{\gamma}&=&2c^{-4}\,dU\otimes dU\\\label{conf2}
\triangle_{\gamma}U&=&0,
\end{eqnarray}
where $\Ric_{\gamma}$ denotes the Ricci curvature tensor of $\gamma$ and $\triangle_{\gamma}$ denotes the $\gamma$-covariant Laplacian on $M^3$.

$N$ was assumed to converge to $1$ asymptotically at spacelike infinity, so $U$ must asymptotically tend to $0$. Indeed, Kennefick and O'Murchadha \cite{KM} showed that the asymptotic flatness assumptions incorporated into the above definition of a geometrostatic system induce the decay conditions
\begin{eqnarray}\label{gamfall}
\gamma_{ij}&=&\delta_{ij}+\mathcal{O}(r^{-2})\\\label{Ufall}
U&=&-\frac{mG}{r}+\mathcal{O}(r^{-2})
\end{eqnarray}
as $r\to\infty$ in suitable coordinates. Here, $m$ is the ADM-mass of the slice $(M^3,g)$, see \cite{ADM,Bartnik,Chrusciel}. In \cite{Cederbaum}, we prove asymptotic estimates in weighted Sobolev spaces that improve this fall-off result. In particular, we find that
\begin{eqnarray}\label{mygamfall}
\gamma_{ij}&=&\left(1-\frac{M^2}{r^2}\right)\delta_{ij}+\frac{2M^2x_{i}x_{j}}{r^4}+\mathcal{O}(r^{-3})\\\label{myUfall}
U&=&-\frac{mG}{r}-\frac{mG\vec{z}_{A}\cdot\vec{x}}{r^3}+\mathcal{O}(r^{-3}),
\end{eqnarray}
as $r\to\infty$, where $M=mG/c^2$ and $\vec{z}_{A}\in\mathbb{R}^3$ is a fixed vector. This decay occurs in asymptotically flat $\gamma$-harmonic coordinates. As a matter of fact, these coordinates coincide with the asymptotically flat (spatial) wave-harmonic coordinates on $(\mathbb{R}\times M^3,ds^2=-N^2c^2dt^2+g)$. The vector $\vec{z}_{A}$ can be interpreted as the coordinate vector of the \emph{asymptotic center of mass} of the system, see below.

In terms of the pseudo-Newtonian variables $\gamma$ and $U$ and inspired by Newtonian gravity, we suggest the following quasi-local definitions of \emph{pseudo-Newtonian mass and center of mass}\footnote{We note that our discussion of center of mass only applies to systems with non-vanishing mass.} of a geometrostatic system $(M^3,g,N)$ with associated pseudo-Newtonian variables $(\gamma,U)$
\begin{eqnarray}\label{m}
 m_{PN}(\Sigma)&:=&\frac{1}{4\pi G}\int_\Sigma\frac{\partial U}{\partial\nu}\,d\sigma\\\label{CoM}
\vec{z}_{PN}(\Sigma)&:=&\frac{1}{4\pi Gm}\int_\Sigma\left(\frac{\partial
U}{\partial\nu}\vec{x}-U\frac{\partial\vec{x}}{\partial\nu} \right)d\sigma,
\end{eqnarray}
where $\Sigma$ is any surface enclosing the support of the matter, $\nu$ and
$d\sigma$ are the outer unit normal to and area measure of $\Sigma$ with respect
to $\gamma$, and $\vec{x}$ is the vector of $\gamma$-harmonic coordinates.

Surprisingly, both of these expressions are \emph{independent} of the particular choice of surface $\Sigma$ (as long as the surface encloses the support of the matter, imagine for example a large coordinate sphere). For the pseudo-Newtonian mass \eqref{m}, this independence of the surface can be seen by combining Equation \eqref{conf2} with the divergence theorem. For the center of mass \eqref{CoM}, the independence of the surface follows from Equation \eqref{conf2} combined with Green's formula and the fact that the coordinates are $\gamma$-harmonic such that $\triangle_{\gamma}\vec{x}=\vec{0}.$ We will thus drop the explicit reference to the surface when referring to pseudo-Newtonian mass and center of mass.

Using the asymptotic decay \eqref{gamfall}, \eqref{Ufall}, we find that
\begin{equation}
m_{PN}=m_{ADM}.
\end{equation}
The total mass of a geometrostatic system is thus \emph{localized}. It can be read off \emph{exactly} on any surface enclosing the matter. Applying Formula \eqref{m} to any smooth surface $\Sigma\subset M^3$, we immediately obtain a notion of mass for an arbitrary part of the system (namely the part bounded by the surface $\Sigma$). By the divergence theorem and \eqref{conf2}, the masses of all components in a multi-component system add up to the total mass of the system just as Newtonian masses do.

If we combine the asymptotic decay \eqref{mygamfall}, \eqref{myUfall} with Formula \eqref{CoM} defining the pseudo-Newtonian center of mass, we find that $\vec{z}_{PN}=\vec{z}_{A}.$ We claim that this vector can indeed be physically interpreted as the coordinate vector of the total center of mass of the system. For this, we exploit Huang's work \cite{Huang} showing that the ADM-center of mass \cite{ADM} coincides\footnote{under precise fall-off conditions at spacelike infinity that are satisfied here.} with the CMC-center of mass constructed via a constant mean curvature (CMC) foliation near infinity by Huisken and Yau \cite{Huisken-Yau} and generalized by Metzger \cite{Metzger}. Using again \eqref{mygamfall} and \eqref{myUfall}, we obtain
\begin{equation}
\vec{z}_{PN}=\vec{z}_{A}=\vec{z}_{ADM}=\vec{z}_{CMC}
\end{equation}
which justifies the name center of mass for the quantities $\vec{z}_{PN}$ and $\vec{z}_{A}$. Moreover, this shows that the center of mass of a geometrostatic system is also \emph{localized}. As above, we obtain a notion of center of mass for an arbitrary part of the system (namely the part bounded by the surface $\Sigma$). By Green's formula, $\gamma$-harmonicity of the coordinates, and \eqref{conf2}, the centers of mass of all components in a multi-component system add up to the total center of mass of the system just as Newtonian centers of mass do.

\section{The Newtonian Limit of Geometrostatics}\label{limit}
Intuitively, the mass and center of mass of a relativistic system should converge to the Newtonian mass and center of mass of its Newtonian limit $c\to\infty$. To make this precise, we use Ehlers' frame theory \cite{Ehlers} which unifies general relativity (GR) and Newton-Cartan gravity (NC) into a common geometric framework with geometric variables $g,h,\Gamma$ and matter tensor $T$.

\begin{figure}[h]
\begin{center}
\resizebox{30ex}{!}{\input{universe.pdftex_t}}\quad\quad\quad\resizebox{30ex}{!}{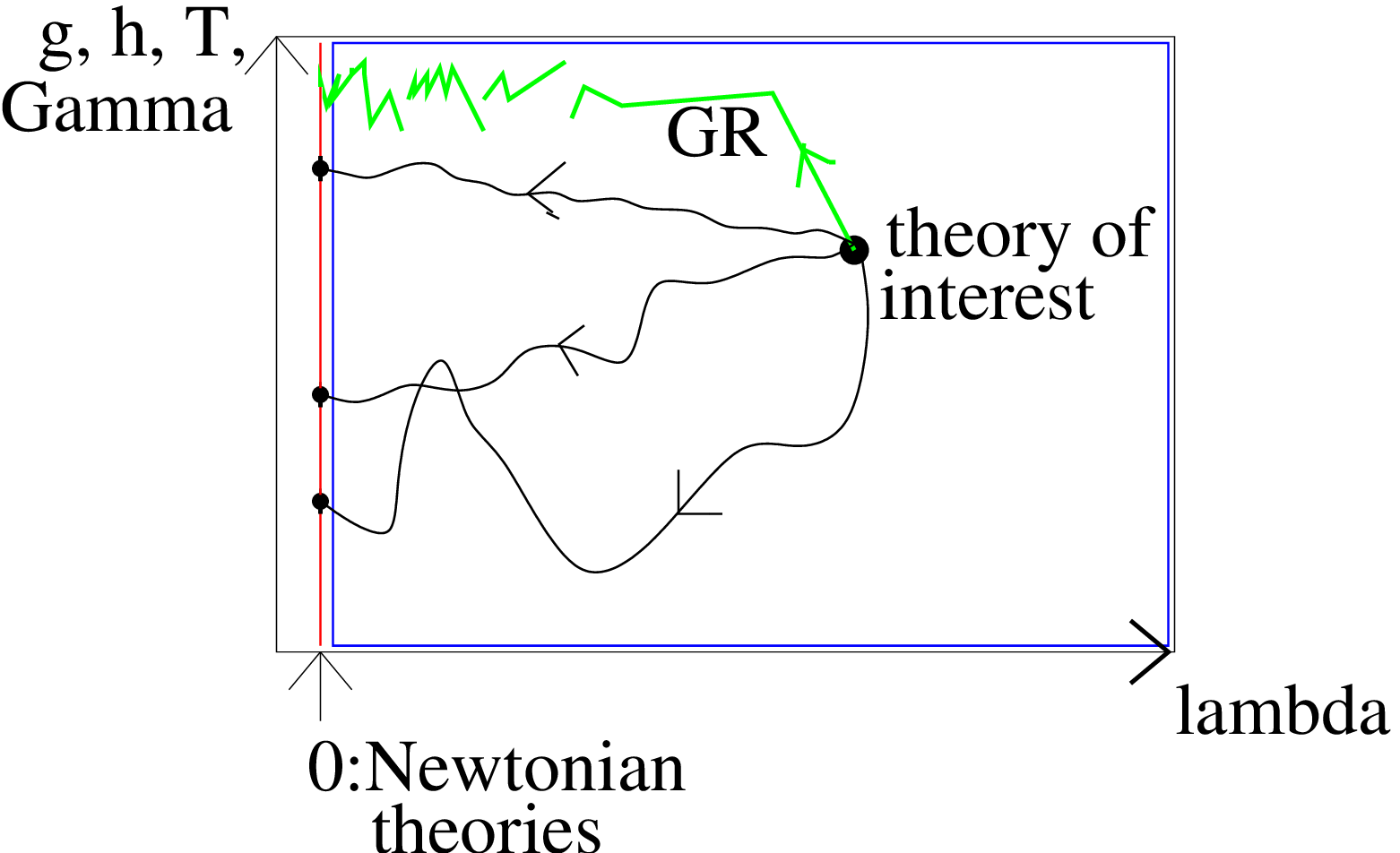}
\caption{\label{limitpic}The universe of Ehlers' frame theory and the Newtonian limit.}
\end{center}
\end{figure}

In frame theory, taking the Newtonian limit corresponds to taking a parametric curve of solutions of GR with parameter $\lambda=c^{-2}$ to its limit $\lambda\to0$, see Figure \ref{limitpic}. Modeling Killing vector fields, staticity, and asymptotic flatness in frame theory, we show that the pseduo-Newtonian potential converges to the Newtonian potential and the metric $\gamma$ converges to the flat metric along any family of geometrostatic systems that possesses a static Newtonian limit. As the localized pseudo-Newtonian formulas \eqref{m} and \eqref{CoM} are nearly identical with the Newtonian formulas, this proves that indeed the relativistic mass and center of mass converge to their Newtonian counterparts.

\section{Conclusion and Outlook}
Geometrostatic systems share many features with Newtonian ones. First, the level sets of the lapse function $N$ (or, equivalently, those of the pseudo-Newtonian potential $U=c^2\ln N$) have the same equipotential properties as the level sets of the Newtonian potential. We thus define the \emph{force on a unit mass test particle} as $\vec{F}:=-\nabla_{\gamma}U$ where $\nabla_{\gamma}$ denotes the $\gamma$-covariant gradient. A second Newtonian type law holds for this notion of force \cite{Cederbaum}.

Secondly, the total mass and center of mass of a geometrostatic system are localized. We put forward explicit geometric formulas for them that also allow for the computation of the notions of mass and center of mass of individual regions. We applied these formulas to prove consistence of ADM-mass and ADM/CMC-center of mass with the Newtonian limit.

This fact and the notion of force might turn out useful for the study of the well-known static $n$-body problem and of Bartnik's conjecture \cite{Szabados}.

\section*{References}
\bibliography{ae100prg}

\end{document}